\documentstyle[12pt]{article}
\makeatletter

\def\newremark#1{\@ifnextchar[{\@ormkwsq{#1}}{\@nrmkwsq{#1}}}

\def\@nrmkwsq#1#2{%
\@ifnextchar[{\@xnrmkwsq{#1}{#2}}{\@ynrmkwsq{#1}{#2}}}

\def\@xnrmkwsq#1#2[#3]{\expandafter\@ifdefinable\csname #1\endcsname
{\@definecounter{#1}\@addtoreset{#1}{#3}%
\expandafter\xdef\csname the#1\endcsname{\expandafter\noexpand
  \csname the#3\endcsname \@rmkcountersep \@rmkcounter{#1}}%
\global\@namedef{#1}{\@rmkwsq{#1}{#2}}\global\@namedef{end#1}
{\@endremarkwithsquare}}}

\def\@ynrmkwsq#1#2{\expandafter\@ifdefinable\csname #1\endcsname
{\@definecounter{#1}%
\expandafter\xdef\csname the#1\endcsname{\@rmkcounter{#1}}%
\global\@namedef{#1}{\@rmkwsq{#1}{#2}}\global\@namedef{end#1}
{\@endremarkwithsquare}}}

\def\@ormkwsq#1[#2]#3{\expandafter\@ifdefinable\csname #1\endcsname
  {\global\@namedef{the#1}{\@nameuse{the#2}}%
\global\@namedef{#1}{\@rmkwsq{#2}{#3}}%
\global\@namedef{end#1}{\@endremarkwithsquare}}}

\def\@rmkwsq#1#2{\refstepcounter
    {#1}\@ifnextchar[{\@yrmkwsq{#1}{#2}}{\@xrmkwsq{#1}{#2}}}

\def\@xrmkwsq#1#2{\@beginremark{#2}{\csname the#1\endcsname}\ignorespaces}
\def\@yrmkwsq#1#2[#3]{\@opargbeginremark{#2}{\csname
       the#1\endcsname}{#3}\ignorespaces}

\def\@rmkcounter#1{\noexpand\arabic{#1}}
\def\@rmkcountersep{.}
\def\@beginremark#1#2{\trivlist \item[\hskip \labelsep{\bf #1\ #2.}]}
\def\@opargbeginremark#1#2#3{\trivlist
      \item[\hskip \labelsep{\bf #1\ #2\ (#3)}]}
\def\@endremarkwithsquare{~\hspace{\fill}~$\Box$\endtrivlist}

\def\@eqnnum{\hbox to .01pt{}\rlap{\rm \hskip -\displaywidth(\theequation)}}
\makeatother

\makeatletter
   \def\@begintheorem#1#2{\sl \trivlist \item[\hskip \labelsep{\bf #2\ #1}]}
   \def\@opargbegintheorem#1#2#3{\sl \trivlist
            \item[\hskip \labelsep{\bf #2\ #1\ (#3)}]}

   \def\section{\@startsection {section}{1}{\z@}{-3.5ex plus -1ex minus
    -.2ex}{2.3ex plus .2ex}{\large\bf}}

   \def\subsection{\@startsection{subsection}{2}{\z@}{-3.25ex plus -1ex minus
   -.2ex}{1.5ex plus .2ex}{\normalsize\bf}}
\makeatother

\makeatletter
\newenvironment{eqn}{\refstepcounter{subsection}
$$}{\leqno{\rm(\thesubsection)}$$\global\@ignoretrue}
\makeatother

\makeatletter
\newenvironment{subeqn}{\refstepcounter{subsubsection}
$$}{\leqno{\rm(\thesubsubsection)}$$\global\@ignoretrue}
\makeatother



\makeatletter
\def\@rmkcounter#1{\noexpand\arabic{#1}}
\def\@rmkcountersep{.}
\def\@beginremark#1#2{\trivlist \item[\hskip \labelsep{\bf #2\ #1.}]}
\def\@opargbeginremark#1#2#3{\trivlist
      \item[\hskip \labelsep{\bf #2\ #1\ (#3).}]}
\def\@endremarkwithsquare{~\hspace{\fill}~$\Box$\endtrivlist}
\makeatother

\newenvironment{prf}[1]{\trivlist
\item[\hskip \labelsep{\it
#1.\hspace*{.3em}}]}{~\hspace{\fill}~$\Box$\endtrivlist}

\newenvironment{proof}{\begin{prf}{\bf Proof}}{\end{prf}}

\let\tempcirc=\circ
\def\circ{\mathord{\raise0.25ex\hbox{$\scriptscriptstyle\tempcirc$}}}

\newcommand{\ZZ}{{\bf Z}}
\newcommand{\QQ}{{\bf Q}}
\newcommand{\RR}{{\bf R}}
\newcommand{\CC}{{\bf C}}
\newcommand{\FF}{{\bf F}}

\newcommand{\TT}{{\bf T}}

\newcommand{\GL}{{\rm GL}}
\newcommand{\Qbar}{\overline{\QQ}}
\newcommand{\End}{{\rm End}}
\newcommand{\Gal}{{\rm Gal}}
\newcommand{\rH}{{\rm H}}

\newcommand{\crys}{{\rm crys}}
\newcommand{\eps}{\varepsilon}
\newcommand{\Frob}{{\rm Frob}}
\newcommand{\Fil}{{\rm Fil}}
\newcommand{\DR}{{\rm DR}}
\newcommand{\ld}{\langle}
\newcommand{\rd}{\rangle}
\newcommand{\Dcrys}{{{\rm D}_\crys}}
\newcommand{\new}{{\rm new}}
\newcommand{\old}{{\rm old}}
\newcommand{\discr}{{\rm discr}}
\newcommand{\trace}{{\rm trace}}
\newcommand{\Ar}{{\rm Ar}}
\newcommand{\vol}{{\rm vol}}
\newcommand{\Spec}{{\rm Spec}}

\newtheorem{theorem}[subsection]{Theorem.}

\newtheorem{corollary}[subsection]{Corollary.}

\newremark{remark}[subsection]{Remark}
\newremark{subremark}[subsubsection]{Remark}
\newremark{notation}[subsection]{Notation}
\newremark{example}[subsection]{Example}
\newremark{definition}[subsection]{Definition}

\begin{document}
\title{On the semi-simplicity of the $U_p$-operator on modular forms.}
\author{Robert F.\ Coleman \and Bas Edixhoven\thanks{partially supported by
the Institut Universitaire de France}}
\maketitle

\section{Introduction.}\label{section1}
For $N$ and $k$ positive integers, let $M^0(N,k)_\CC$ denote the $\CC$-vector
space of cuspidal modular forms of level $N$ and weight~$k$. This vector
space is equipped with the usual Hecke operators $T_n$, $n\geq1$. If we
need to consider several levels or weights at the same time, we will denote
this $T_n$ by $T_n^N$, or~$T_n^{N,k}$. If $p$ is a prime number dividing $N$,
our $T_p$ is also known under the name~$U_p$. One of our main results can
be stated very easily: if $k=2$ and $p^3$ does not divide $N$, then the
operator $T_p$ is semi-simple. We can prove the same result for weight
$k\geq3$, under the assumption that certain crystalline Frobenius
elements are semi-simple. Milne has shown in \cite[\S2]{Milne1} that this
semi-simplicity is implied by Tate's conjecture claiming
that for $X$ projective and smooth over a finite field of characteristic $p$,
and $r\geq0$, $\dim_\QQ({\rm CH}^r(X)/{\equiv}_{\rm num})$ equals the order
of $\zeta(X,s)$ at~$r$.
Ulmer proved in \cite{Ulmer1} that $T_p$ is semi-simple, for $k=3$ and
$p^2$ not dividing~$N$, under the assumption of the Birch-Swinnerton-Dyer
conjecture for elliptic curves over function fields in characteristic~$p$.
His method is quite different from ours: assuming that $T_p$ is not
semi-simple, he really shows that the Birch-Swinnerton-Dyer conjecture
does not hold for an explicitly given elliptic curve.

The structure of our proof is as follows. Using the theory of newforms, the
problem is shown to be equivalent to the problem of showing that, for a
normalized
newform $f$ of weight $k$, prime-to-$p$ level and character $\eps$, the
polynomial $x^2-a_px+\eps(p)p^{k-1}$ has no double root. This polynomial
happens to be the characteristic polynomial of the Frobenius element at $p$
in the two-dimensional Galois representations associated to~$f$; it is also
the characteristic polynomial of the crystalline Frobenius asociated to~$f$.
We show that this crystalline Frobenius cannot be a scalar.

In Sections~\ref{section2} and \ref{section3} we prove the results concerning
these Frobenius elements for $k=2$ and $k\geq2$, respectively.
Section~\ref{section2} is quite elementary, whereas in Section~\ref{section3}
we use a lot of the machinary for comparing $p$-adic etale and crystalline
cohomology.
In Section~\ref{section4} we give some applications: the Ramanujan inequality
is a strict inequality in certain cases, certain Hecke algebras are reduced,
hence have non-zero discriminant. Section~\ref{section5} gives some
results, due to Abbes and Ullmo, concerning the discriminants of certain
Hecke algebras.

To end this introduction, let us explain why the case $k=1$ is completely
different. Consider a normalized cuspidal eigenform $f=\sum a_nq^n$ of some
level $N$, of weight one and with some character~$\eps$. Deligne and Serre
have shown (\cite[\S4]{DeligneSerre})
that there exists a continuous representation
$\rho_f$ from $\Gal(\Qbar/\QQ)$ to $\GL_2(\CC)$, unramified outside $N$, such
that, for all primes $p$ not dividing $N$, the characteristic polynomial of a
Frobenius element at $p$ is $x^2-a_px+\eps(p)$. Since the image of $\rho_f$
is finite, Chebotarev's density theorem gives the existence of primes $p$
not dividing $N$ such that the Frobenius element at $p$ is the identity
element, hence has characteristic polynomial $(x-1)^2$.

\section{An elementary proof in the case of weight two.}\label{section2}

\begin{theorem}\label{thm2.1}
Let $f=\sum a_nq^n$ be a cuspidal normalized eigenform of weight two, some
level $N$ and character $\eps\colon(\ZZ/N\ZZ)^*\to\CC^*$. Let $p$ be a
prime number not dividing~$N$. Then the polynomial $x^2-a_px+\eps(p)p$ has
simple roots.
\end{theorem}
\begin{proof}
The proof is by contradiction, so we suppose that the polynomial has a
double root~$\lambda$. Then of course we have $\lambda^2=\eps(p)p$ and
$2\lambda=a_p$. Let $K$ be the finite extension of $\QQ$ generated by
the $a_n$ and the $\eps(a)$, and let $O_K$ be its ring of integers.
Let $J$ denote the
jacobian variety of the modular curve $X_1(N)$ over~$\QQ$. We identify
the space $M^0(N,2)_\QQ$ of weight two cuspforms of level $N$ and with
coefficients in $\QQ$ with the cotangent space at the origin of~$J$; this is
compatible  with the action of the Hecke operators. Let $\TT$ be the
subring of $\End(J)$ that is generated by the $T_n$, $n\geq1$, and the
diamond operators $\ld a\rd$, $a\in(\ZZ/N\ZZ)^*$. Let $I$ be
the annihilator of $f$ in $\TT$, and let $A'_\QQ:=J/IJ$ be the quotient of
$J$ by its subvariety generated by the images of all elements of~$I$.
It is well known that $A'_\QQ$ has dimension $[K:\QQ]$, and that for
every prime number $l$, the free $K\otimes\QQ_l$-module $V_l(A'_\QQ)$ of rank
two gives the $l$-adic Galois representation $\rho_{f,l}$ associated to~$f$.
We prefer to work with an abelian variety $A_\QQ$ that is isogeneous to
$A'_\QQ$ and on which we have an action of all of~$O_K$. This is easily done:
define $A_\QQ:=O_K\otimes_\TT A'_\QQ$, where the tensor product should be
calculated by taking a presentation of~$O_K$.

The abelian variety $A_\QQ$ has good reduction at $p$; let $A_{\ZZ_p}$
denote the corresponding abelian scheme over~$\ZZ_p$. Consider the
first algebraic de Rham cohomology group $M:=\rH^1_\DR(A_{\ZZ_p}/\ZZ_p)$.
It is a free $\ZZ_p$-module of rank $2[K:\QQ]$, equipped with its Hodge
filtration:
\begin{eqn}\label{eqn2.2}
M = \Fil^0M \supset \Fil^1 M = \rH^0(A_{\ZZ_p},\Omega^1).
\end{eqn}
The submodule $\Fil^1M$ is free of rank $[K:\QQ]$ as $\ZZ_p$-module, and
has the property that $\Fil^0M/\Fil^1M$ is torsion free. The double root
$\lambda$ of $x^2-a_px+\eps(p)p$ is in $O_K$, since it is integral and
$2\lambda$ is in~$K$. In the endomorphism ring of $A_{\FF_p}$ we have the
Eichler-Shimura congruence relation:
\begin{eqn}\label{eqn2.3}
0 = (\Frob_p-\Frob_p)(\Frob_p-\Frob_p') = \Frob_p^2-a_p\Frob_p+\eps(p)p
= (\Frob_p-\lambda)^2,
\end{eqn}
where $\Frob_p$ is the Frobenius endomorphism and $\Frob_p'$ the
Verschiebung, multiplied by~$\eps(p)$. The fact that every abelian variety
over $\FF_p$ is isogeneous to a product of simple ones implies that $\Frob_p$
is semi-simple in the sense that it satisfies an identity of the form
$P(\Frob_p)=0$, with $P$ a polynomial with coefficients in $\QQ$ that has
simple roots. It follows that $\Frob_p=\lambda$ in $\End(A_{\FF_p})$.
Since $O_K\otimes\ZZ_p$ is a product of a finite number of discrete
valuation rings, $\Fil^1M$ is a locally free module over it; it is in
fact free of rank one (consider $\QQ\otimes\Fil^1M$). It follows that
$\lambda$ does not annihilate $\FF_p\otimes\Fil^1M$, since we have
$\lambda^2=\eps(p)p$. But $\FF_p\otimes\Fil^1M$ is the same as
$\rH^0(A_{\FF_p},\Omega^1)$, and on this module $\lambda$ acts as
$\Frob_p^*$, hence it does annihilate. This contradiction finishes the proof.
\end{proof}

\section{The general case.}\label{section3}
In this section we try to generalize Theorem~\ref{thm2.1} as much as we can
to higher weights. For doing that we replace the module $M$ of
Section~\ref{section2} by the $p$-adic crystalline realization of the motive
associated to~$f$; this gives us a filtered $\phi$-module $M$ of rank two.
The comparison theorem for crystalline and $p$-adic etale cohomology implies
that this filtered $\phi$-module is weakly admissible, from which it follows
immediately that the crystalline Frobenius $\phi$ cannot be a scalar.
Unfortunately, it is not known that $\phi$ is semi-simple, so all we show is
that semi-simplicity of $\phi$ implies that the polynomial
$x^2-a_px+\eps(p)p^{k-1}$ has simple roots.

Let $f=\sum a_nq^n$ be a normalized cuspidal newform of some level $N$,
weight $k\geq2$ and character~$\eps$.
Let $K$ be the field generated by the $a_n$ and the~$\eps(a)$.
Let $p$ be a prime number not dividing~$N$. Our first objective
is to construct the $p$-adic crystalline realization of the motive associated
to~$f$. In \cite{Scholl1}, Scholl constructs a Grothendieck motive $M(f)$
over $\QQ$, with coefficients in $K$, such that for every prime number $l$
the Galois representation $\rho_{f,l}\colon G_\QQ\to\GL_2(\QQ_l\otimes K)$
is the dual of the $l$-adic realization~$\rH_l(M(f))$. Concretely, he
constructs a projector in the group ring of a finite group of automorphisms
of the smooth and projective model $X$ of the $k{-}2$-fold fibered product of
the universal elliptic curve over $Y(N')$ (with $N'$ a suitable multiple of
$N$) constructed by Deligne in \cite[\S5]{Deligne1},
such that the $l$-adic and Betti realizations of the corresponding Chow motive
are, in a way that is compatible with Hecke operators, the parabolic
cohomology groups used in~\cite{Deligne1}. The Grothendieck motive associated
to this Chow motive (i.e., one replaces rational equivalence by homological
equivalence) has an action by the Hecke algebra of the space of cuspidal
modular forms of weight $k$ on the modular curve $X(N')_\QQ$; $M(f)$ is a
suitable factor. The variety $X$ over $\QQ$ has a smooth projective
model over $\ZZ[1/N']$ (see \cite[4.2.1]{Scholl1} and \cite{Deligne1}),
hence $M(f)$ has a crystalline realization
$M:=\rH_\crys(M(f))$ which is a free $\QQ_p\otimes K$-module of rank two
equipped with an endomorphism $\phi$, the crystalline Frobenius, that is
induced by the Frobenius endomorphism of the reduction mod $p$ of~$X$.
The characteristic polynomial of $\phi$ is $x^2-a_px+\eps(p)p^{k-1}$;
this can be shown in the same way as one can show it for the $l$-adic
realizations, or one invokes a result of Katz and Messing
(see \cite[4.2.3]{Scholl1}).

\begin{theorem}\label{thm3.1}
Let $f=\sum a_nq^n$ be a normalized cuspidal newform of some level $N$,
weight $k\geq2$ and character~$\eps$. Let $p$ be a prime number not
dividing~$N$. Then the Frobenius $\phi$ of the crystalline realization $M$
of the motive $M(f)$ is not scalar, i.e., it is not in $\QQ_p\otimes K$.
\end{theorem}
\begin{proof}
The proof is by contradiction. We suppose that $\phi$ is an element,
$\lambda$ say, of $\QQ_p\otimes K$. The comparison theorem for crystalline
and de Rham cohomology for smooth proper $\ZZ_p$-schemes
(see \cite[\S1.3]{Illusie1}) gives us an
isomorphism between $M$ and $\QQ_p\otimes \rH_\DR(M(f))$, and hence a
decreasing filtration (the Hodge filtration) $\Fil$ on~$M$. So
$(M,\phi,\Fil)$ is an object of the category of filtered $\phi$-modules.
(A filtered $\phi$-module is a finite dimensional $\QQ_p$-vector space $M$
with a decreasing, exhaustive and separating filtration $\Fil^i$, $i\in\ZZ$,
and an endomorphism $\phi$; morphisms are linear
maps respecting $\Fil$ and~$\phi$; see~\cite[\S2.3]{Illusie1}.)
The Hodge filtration on $\rH_\DR(M(f))$ induces the Hodge decomposition
of $\CC\otimes\rH_{\rm B}(M(f))$, which is of type $(k{-}1,0),(0,k{-}1)$,
hence $\Fil^0(M)=M$,
$\Fil^1(M)=\Fil^{k-1}(M)$ is free of rank one, and $\Fil^k(M)=0$.
Fontaine has constructed Grothendieck's ``mysterious functor''
$\Dcrys$ from the category of finite dimension $\QQ_p$-vector spaces with
continuous $\Gal(\Qbar_p/\QQ_p)$-action to the category of filtered
$\phi$-modules (see the introduction of~\cite{Illusie1}). It is a theorem
of Faltings (see \cite[\S3.2]{Illusie1}), of which a special case was proved
earlier by Fontaine and Messing, that there is an isomorphism
of filtered $\phi$-modules between $M$ and $\Dcrys(\rH_p(M(f)))$. In fact,
the theorem is stated for varieties, but since the isomorphism is compatible
with the multiplicative structure and with cycle classes, it also works for
Grothendieck motives.

The most important consequence of this theorem for us
is that the filtered $\phi$-module $M$ is admissible, hence weakly
admissible, in the sense of Fontaine, see \cite[4.4.6]{Fontaine1}.
Recall that to a filtered $\phi$-module $M$ one associates two polygons:
the Hodge polygon, depending only on the filtration, and the Newton polygon,
depending only on~$\phi$. Weakly admissible means that for every subobject
$M'$ of $M$ the Newton polygon lies above or on the Hodge polygon, and that
the two polygons for $M$ itself have the same end-point. An equivalent
formulation is the following. For $M$ a filtered $\phi$-module let
$t_N(M)$ be the $p$-adic valuation of the determinant of $\phi$, and let
$t_H(M)$ be the maximal $i$ such that $\Fil^i(\det M)\neq0$. Then $M$ is
weakly admissible if and only if firstly: $t_N(M)=t_H(M)$, and
secondly: for all subobjects $M'$ of $M$ one has $t_H(M')\leq t_N(M')$.

Consider now our weakly admissible filtered $\phi$-module~$M$.
Since $\phi$ is the element $\lambda$ of $\QQ_p\otimes K$, we have the
subobject $M':=\Fil^{k-1}(M)$ of $M$ (we give it the induced filtration).
Then $t_H(M')=[K:\QQ](k-1)$, whereas $t_N(M')=[K:\QQ](k-1)/2$ (recall that
$M'$ is free of rank one over $\QQ_p\otimes K$ and that
$\lambda^2=\eps(p)p^{k-1}$). Since $k\geq2$, this contradicts the weak
admissibility of~$M$.
\end{proof}

\begin{corollary}\label{cor3.2}
Let $f=\sum a_nq^n$ be a normalized cuspidal eigenform of some level $N$,
weight $k\geq2$ and character~$\eps$. Let $p$ be a prime number not
dividing $N$ and suppose that the crystalline Frobenius $\phi$ of the
$\QQ_p$-vector space $M(f)$ is semi-simple. Then the polynomial
$x^2-a_px+\eps(p)p^{k-1}$ has simple roots.
\end{corollary}
\begin{proof}
The proof is by contradiction: we suppose that $\lambda\in K$ is a double
root. As we have already said above, the polynomial in question is the
characteristic polynomial of the endomorphism $\phi$ of the free rank two
$\QQ_p\otimes K$-module~$M$. Hence it satisfies the identity
$(\phi-\lambda)^2=0$. Now $\QQ_p\otimes K$ is a finite product of fields,
hence the semi-simplicity of $\phi$ implies that $\phi$ is multiplication
by~$\lambda$. But this contradicts Theorem~\ref{thm3.1}.
\end{proof}

\begin{remark}\label{remark3.3}
The first three lines of \cite[\S2]{Milne1} show that Tate's conjecture
mentioned in Section~\ref{section1} implies the semi-simplicity of
$l$-adic and crystalline Frobenius elements of smooth projective
varieties over finite fields.
\end{remark}

\begin{remark}\label{remark3.4}
Scholl remarks that his explicit construction of the crystalline realization
$M$ of $M(f)$ in \cite{Scholl2} should show directly that $M$ is weakly
admissible.
\end{remark}

\section{Applications.}\label{section4}
\begin{theorem}\label{thm4.0}
Let $N\geq1$ and $k\geq2$ be integers. Let $f=\sum a_nq^n$ be a normalized
cuspidal eigenform of level $N$ and weight~$k$. Let $p$ be a prime number
not dividing~$N$. If $k>2$ assume Tate's conjecture mentioned in
Section~\ref{section1}. Then we have $|a_p|<2p^{(k-1)/2}$.
\end{theorem}
\begin{proof}
Let $\eps$ be the character of~$f$. Theorems~\ref{thm2.1} and~\ref{thm3.1}
show that $x^2-a_px+\eps(p)p^{k-1}$ has no double root. Deligne has shown
(\cite[Thm.~6.1]{DeligneSerre} and \cite[Thm.~1.6]{Deligne2})
that the roots have absolute value $p^{(k-1)/2}$.
\end{proof}

\begin{theorem}\label{thm4.1}
Let $N\geq1$ and $k\geq2$ be integers. Let $p$ be a prime number such that
$p^3$ does not divide~$N$. Assume Tate's conjecture mentioned in
Section~\ref{section1} if $k\geq3$ and $p|N$. Then the endomorphism $T_p$ of
$M^0(N,k)_\CC$ is semi-simple.
\end{theorem}
\begin{proof}
For the sake of notation, let $N\ge1$ be an integer and $p$ a prime number
not dividing~$N$. Let $k\geq1$. In this case $T_p$ is normal with respect to
the Petersson scalar product on $M^0(N,k)_\CC$. Hence $T_p$ is diagonalizable.

Let us now consider $M^0(pN,k)_\CC$, with $k\geq2$. By the theory of
newforms,
$M^0(pN,k)_\CC$ is the direct sum of its $p$-new part $M^0(pN,k)_\CC^{p\new}$
and its $p$-old part $M^0(pN,k)_\CC^{p\old}$, this decomposition being
respected by all Hecke operators. The restriction of $T^{pN,k}_p$ to
$M^0(pN,k)_\CC^{p\new}$ is normal, hence diagonalizable. The
$p$-old part is isomorphic to the direct sum of two copies of
$M^0(N,k)_\CC$, via the map $(f(q),g(q))\mapsto f(q)+g(q^p)$. The restriction
of $T^{pN,k}_p$ to $M^0(pN,k)_\CC^{p\old}$ is then given by the following
two by two matrix:
\begin{subeqn}\label{eqn4.1.1}
T^{pN,k}_p|_{M^0(pN,k)_\CC^{p\old}} =
\left(
\begin{array}{cc}
T^{N,k}_p & 1\\ -p^{k-1}\ld p\rd & 0
\end{array}
\right).
\end{subeqn}
We have already seen that $M^0(N,k)_\CC$ is the direct sum of its common
eigenspaces $V_{a_p,\eps}$ for $T_p$ and the diamond operators. It follows
that $M^0(pN,k)_\CC^{p\old}$ decomposes as a direct sum of terms
$V_{a_p,\eps}^2$, and that the restriction of $T^{pN,k}_p$ to each of the
$V_{a_p,\eps}^2$ is annihilated by $x^2-a_px+\eps(p)p^{k-1}$. Under the
hypotheses of the theorem we are proving, these polynomials have simple
roots by Theorems~\ref{thm2.1} and~\ref{thm3.1}.

Let us now consider $M^0(p^2N,k)_\CC$, with $k\geq2$. Here too this space
is the direct sum of its $p$-old and $p$-new parts. On the $p$-new part
$T^{p^2N,k}_p$ is self-adjoint, hence diagonalizable. The $p$-old part is
now isomorphic to the direct sum of three copies of $M^0(N,k)_\CC$
and two copies of $M^0(pN,k)_\CC^{p\new}$.
The restrictions of $T^{p^2N,k}_p$ to $M^0(N,k)_\CC^3$ and
$(M^0(pN,k)_\CC^{p\new})^2$ are given by the following matrices:
\begin{subeqn}\label{eqn4.1.2}
\left(
\begin{array}{ccc}
T^{N,k}_p & 1 & 0 \\
-p^{k-1}\ld p\rd & 0 & 1\\
0 & 0 & 0
\end{array}
\right), \qquad
\left(
\begin{array}{cc}
T^{pN,k}_p & 1\\
0 & 0
\end{array}
\right)
\end{subeqn}
One can now repeat the same type of argument as above, invoking
Theorems~\ref{thm2.1} and \ref{thm3.1} to see that, under the hypotheses
of the theorem we are proving, $x(x^2-a_px+\eps(p)p^{k-1})$, with
$(a_p,\eps(p))$ as before, has simple roots. The space
$M^0(pN,k)_\CC^{p\new}$ is a direct sum of eigenspaces for $T^{Np,k}_p$,
and one knows that the eigenvalues are non-zero
(see~\cite[\S1.8]{DeligneSerre}).
It follows that also the restriction of $T^{p^2N,k}_p$ to
$(M^0(pN,k)_\CC^{p\new})^2$ is diagonalizable.
\end{proof}

\begin{corollary}\label{cor4.2}
Let $N\geq1$ be cube free, and let $k\geq2$. Let $\TT$ be the $\ZZ$-algebra
generated by the endomorphisms $T_n$, $n\geq1$, and $\ld a\rd$,
$a\in(\ZZ/N\ZZ)^*$, of $M^0(N,k)_\CC$. Assume Tate's conjecture mentioned in
Section~\ref{section1} if $k>2$. Then the ring $\TT$ is reduced.
\end{corollary}
\begin{proof}
This is so because $\TT$ is a subring of the $\CC$-algebra
$\TT_\CC:=\CC\otimes\TT$
generated by the $T_n$ and~$\ld a\rd$. Theorem~\ref{thm4.1} tells us
that the $T_p$ and $\ld a\rd$ can be simultaneously diagonalized.
They generate $\TT_\CC$, hence $\TT_\CC$ is a product of copies of $\CC$,
hence reduced.
\end{proof}

\begin{remark}\label{rmk4.3}
For general $N$ and $k$, the Hecke algebra $\TT$ is well known to be a
free $\ZZ$-module; $\TT_\QQ:=\QQ\otimes\TT$ is well known
to be Gorenstein, i.e., its $\QQ$-linear dual $(\QQ\otimes\TT)^\vee$ is
free of rank one as $\QQ\otimes\TT$-module, see for
example~\cite[p.~481]{Wiles1}.
One way to prove this is as follows. By the $q$-expansion principle,
$M^0(N,k)_\CC^\vee$ is free of rank one as $\TT_\CC$-module. Then one
constructs a $\TT_\CC$-bilinear $\CC$-valued pairing on $M^0(N,k)_\CC$
to get an isomorphism of $\TT_\CC$-modules between $M^0(N,k)_\CC$
and its dual. Another way to prove it is to use the theory of new forms.
This last proof gives more information on how exactly $\TT_\CC$
decomposes as a product of $\CC$-algebras. Parent needed such
information and so he worked out the details in~\cite{Parent1}.
His work made see that
Theorem~\ref{thm4.1} should be stated for $N$ cube free instead of
square free.
Of course, statements that completions of $\TT$ at certain of its
maximal ideals are Gorenstein are much more subtle and harder to
prove (see for example~\cite[\S2.1]{Wiles1}).
\end{remark}

\section{Discriminants of Hecke algebras.}\label{section5}
According to Corollary~\ref{cor4.2}, certain Hecke algebras $\TT$ are
reduced. This means that the discriminants $\discr(\TT)$ of their trace
forms $(x,y)\mapsto \trace(xy)$ are non-zero. These discriminants ``count''
all congruences between different eigenforms of fixed level and weight, hence
are quite useless for dealing with congruences with a fixed form (in
particular, nothing interesting can be said on the degrees of modular
parametrizations of elliptic curves over $\QQ$).
The following result, relating
such discriminants to heights of modular curves, is due to Abbes and Ullmo
(unpublished).
\begin{theorem}[Abbes, Ullmo]\label{thm5.1}
Let $p$ be a prime number, and let $\TT$ be the Hecke algebra associated
to $M^0(p,2)_\CC=\rH^0(X_0(p)_\CC,\Omega)$. Then one has:
$$
h(X_0(p)_\QQ) = \frac{1}{2}\log|\discr(\TT)| - \sum_{i=1}^g \log\|\omega_i\|,
$$
where $h$ is the modular height of curves over $\QQ$
(see \cite[\S3.3]{Szpiro1}), where
$\omega_1,\ldots,\omega_g$ are the normalized eigenforms and $\|\cdot\|$
the norm of the scalar product
$\ld\omega|\eta\rd= (i/2)\int_{X_0(p)(\CC)}\omega\wedge\overline{\eta}$
on $\rH^0(X_0(p)_\CC,\Omega)$.
\end{theorem}
\begin{proof}
We start by recalling the definition of~$h$. So let $X_\QQ$ be a smooth
proper geometrically irreducible curve over~$\QQ$, of some genus~$g$.
Let $J_\QQ$ be its jacobian and $J$ the N\'eron model over~$\ZZ$. Then we
have the free $\ZZ$-module of rank one
$\omega_J:=\Lambda^g0^*\Omega^1_{J/\ZZ}$, with the
scalar product on $\CC\otimes\omega_J$ given by
$\ld\omega|\eta\rd=
(i/2)^g(-1)^{g(g-1)/2}\int_{J(\CC)}\omega\wedge\overline{\eta}$.
The height $h(X_\QQ)$ is then defined to be the Arakelov degree
of this metrized line bundle: $h(X_\QQ)=\deg_\Ar(\omega_J)=-\log\|\omega\|$,
with $\omega$ a generator of~$\omega_J$. The $\ZZ$-module $\omega_J$ is
equipped with the scalar product on
$\CC\otimes\omega_J=\rH^0(X_0(p)_\CC,\Omega)$ already mentioned in the
theorem above. This scalar product induces a real scalar product, and
hence a volume form, on $\RR\otimes\omega_J$ (the volume form being
determined by the condition that a cube with edges of length one has
volume one). A calculation (see \cite[lemme~3.2.1]{Szpiro1}) shows that
one has:
\begin{subeqn}\label{eqn5.1.1}
h(X_\QQ) = -\log\vol(\RR\otimes\omega_J/\omega_J).
\end{subeqn}
Let now $X$ be the curve~$X_0(p)_\QQ$. In that case, $\omega_J$ is the
same as $\rH^0(X_\ZZ,\Omega)$, with $X_\ZZ$ the usual model over $\ZZ$
(semi-stable, $X_{\FF_p}$ consisting of two irreducible components) and
$\Omega$ its dualizing sheaf (\cite[Ch.~II, \S3]{Mazur1}). The fact that the
two irreducible components of $X_{\FF_p}$ are of genus zero implies that
the pairing
\begin{subeqn}\label{eqn5.1.2}
\TT \times \omega_J \to \ZZ,\quad (t,\omega)\mapsto a_1(t\omega),
\end{subeqn}
with $a_1$ denoting the linear form on $\omega_J$ that takes the coefficient
of $q$ in the $q$-expansion, is perfect, i.e., it induces an isomorphism
between $\TT^\vee$ of $\TT$ and $\omega_J$.
Let $\omega_1,\ldots,\omega_g$ be as in the theorem. Sending an element
of $\TT_\RR$ to the eigenvalues of the $\omega_i$ for it is an
isomorphism of $\RR$-algebras:
\begin{subeqn}\label{eqn5.1.3}
\TT_\RR \to \RR^g, \quad t\mapsto (a_1(t\omega_1,\ldots,a_1(t\omega_g)).
\end{subeqn}
The trace form on $\TT$ corresponds to the standard scalar product on~$\RR^g$.
Composing the dual of the isomorphism (\ref{eqn5.1.3}) with the isomorphism
$\TT_\RR^\vee\to\RR\otimes\omega_J$ from (\ref{eqn5.1.2}) gives an
isomorphism from $\RR^g$ to $\RR\otimes\omega_J$ mapping the $i$th standard
basis vector $e_i$ to~$\omega_i$. It follows that the volume form on
$\RR\otimes\omega_J$ corresponds to $\prod_i\|\omega_i\|$ times the one
on $\TT_\RR^\vee$ corresponding to the trace form. We find:
\begin{subeqn}\label{eqn5.1.4}
\vol(\RR\otimes\omega_J/\omega_J) =
\left(\prod_{i=1}^g \|\omega_i\|\right)\vol(\TT_\RR^\vee/\TT)
= \left(\prod_{i=1}^g \|\omega_i\|\right)\vol(\TT_\RR/\TT)^{-1}.
\end{subeqn}
The proof is finished since $\vol(\TT_\RR/\TT)=|\discr(\TT)|^{1/2}$.
\end{proof}

\begin{theorem}[Abbes, Ullmo]\label{thm5.2}
For every $\eps>0$ there exists $c(\eps)$ in $\RR$ such that for all
prime numbers $p$ one has $h(X_0(p)_\QQ)\leq c(\eps)p^{1+\eps}$.
\end{theorem}
\begin{proof}
Let $\TT$ and $g$ be as above, and let $\TT':=\sum_{i=1}^g\ZZ T_i$. Then
$\TT'$ is of finite index in $\TT$ because $\infty$ is not a Weierstrass
point of $X_0(p)$ (see \cite[\S3]{LehnerNewman} or \cite[\S4]{Edixhoven1}).
The image of $T_i$ under the isomorphism~(\ref{eqn5.1.3}) is
$(a_i(\omega_1),\ldots,a_i(\omega_g))$. It follows that we have the
equalities:
\begin{subeqn}\label{eqn5.2.1}
\discr(\TT) = \discr(\TT') |\TT/\TT'|^{-2}, \quad
|\discr(\TT')| = |\det_{1\leq i,j\leq g} a_i(\omega_j)|^2.
\end{subeqn}
Weil's theorem on absolute values of eigenvalues of Frobenius endomorphisms
of abelian varieties over finite fields implies that
$|a_i(\omega_j)|\leq \sigma(i)i^{1/2}$, where $\sigma(i)$ is the number of
positive integers dividing~$i$. It follows that
$|\discr(\TT')|\leq g!\prod_{i=1}^g\sigma(i)i^{1/2}$. The rest of the proof
consists of applying Theorem~\ref{thm5.1} and standard estimates (including
an absolute lower bound for the~$\|\omega_i\|$).
\end{proof}

\begin{remark}\label{rmk5.3}
One knows (see \cite[Ch.~II, Prop.~10.6]{Mazur1})
that $\Spec(\TT)$ (for $\TT$ as above) is connected. This implies a
lower bound for $\discr(\TT)$ (use \cite{Odlyzko1}).
On the other hand, the $\|\omega_i\|$ are bounded above by a constant
times~$p$. Unfortunately, the lower bound for $h(X_0(p)_\QQ)$ obtained
like this seems too weak to be useful. Assuming $\TT$ to be Gorenstein does
not significantly improve this lower bound.
\end{remark}

\begin{remark}\label{rmk5.4}
Several problems arise when one wants to generalize the above results for
$X_0(p)_\QQ$ to more general~$X_0(N)_\QQ$. First of all, $\TT^\vee$ will not
be the same as $\omega_J$, but it should be possible to estimate
$|\TT^\vee/\omega_J|$. Secondly, the comparison of the trace form on $\TT$
and the scalar product on $\RR\otimes\omega_J$ is more complicated. Note that
the trace form can even be degenerate, if $N$ is not cube free. Thirdly,
it seems to be unknown if $\infty$ can be a Weierstrass point on $X_0(N)_\QQ$
when $N$ is square free (see \cite[\S3]{LehnerNewman}; it is known that
$\infty$ is a Weierstrass point when $N$ is divisible by 4 or 9, for example).

Suppose now that $N$ is square free.
Before knowing the result that $\TT$ is reduced, Abbes and Ullmo have related
$h(X_0(N))$ to the discriminants of the ``new parts'' of the Hecke algebras
of level dividing~$N$ (unpublished). The techniques they use come from
\cite{AbbesUllmo1} and~\cite{AbbesUllmo2}. They show that Theorem~\ref{thm5.2}
holds for square free $N$ such that $\infty$ is not a Weierstrass point.
With the same techniques it is certainly possible to solve the first two
of the three problems mentioned above.
\end{remark}

\vspace{1.5cm}\noindent
{\bf Acknowledgements.} We thank Rutger Noot for a discussion
during which we were led to a proof of Theorem~\ref{thm2.1}. Bas Edixhoven
is grateful to Richard Taylor for a discussion concerning the case of
general weights. He also thanks Pierre Berthelot for explaining to him certain
results on Newton and Hodge polygons. Thanks also go to Ahmed Abbes and
Emmanuel Ullmo for communicating to us their results mentioned in
Section~\ref{section5}, and permitting us to publish them.
Bas Edixhoven wants to thank the Centre for Research in Mathematics at the
Institut d'Estudis Catalans in Barcelona where he finally found time to
write this text.
Robert Coleman thanks the University of Rennes for inviting him for a month
in~1995.

\vfill
\begin{center}
\parbox[t]{6cm}
{Robert Coleman\\Dept. of Mathematics\\University of California\\
Berkeley, CA 94720\\U.S.A.}
\hspace{2cm}
\parbox[t]{6cm}
{Bas Edixhoven\\IRMAR\\Campus de Beaulieu\\35042 Rennes cedex\\France}
\end{center}

\end{document}